\documentclass[preprint,12pt]{elsarticle}
\usepackage{amsfonts,amssymb,amscd}
\usepackage{oldgerm}

\usepackage[Conny]{fncychap}
\usepackage[latin1]{inputenc}
\usepackage{amssymb,amsmath,amsthm,amsfonts}
\usepackage[latin1]{inputenc}
\usepackage{mathrsfs} %\mathscr
\usepackage[all]{xy}

\newtheorem{teo}{Theorem}[section]

\newtheorem{obs}{Remark}[section]

\newtheorem{id}[teo]{Identity}
\newtheorem{statement}{Statement}[section]

\newcommand{\C}{{\mathbb C}}
\newcommand{\R}{{\mathbb R}}
\newcommand{\Z}{{\mathbb Z}}
\newcommand{\N}{{\mathbb N}}

\begin{document}
\begin{frontmatter}

\title{Integrating the Jacobian equation}

\begin{abstract}
We show essentially that the differential equation $\frac{\partial (P,Q)}{\partial (x,y)} =c \in {\mathbb C}$, for $P,\,Q \in {\mathbb C}[x,y]$, may be "integrated", in the sense that it is equivalent to an algebraic system of equations involving the homogeneous components of $P$ and $Q$. Furthermore, the first equations in this system give explicitly the homogeneous components of $Q$ in terms of those of $P$. The remaining equations involve only the homogeneous components of $P$.
\end{abstract}


\author[air]{A. S. de Medeiros}
\ead{airtonsoh@yahoo.com.br}

\author[rvt]{R. R. Silva \corref{cor1}}
\ead{rrsilva73@gmail.com}

\cortext[cor1]{Corresponding author. Tel. +55 61 3107-6454}

\address[air]{ Instituto de Matem\'atica, Universidade Federal do Rio de Janeiro, Ilha do Fund\~{a}o, CEP: 21941-909, Rio de Janeiro, RJ, Brazil}

\address[rvt]{Departamento de Matem\'atica, Universidade de Bras\'ilia, Campus Universit\'ario Darcy Ribeiro, Asa Norte, CEP:
70910-900, Bras\'ilia, Brazil}

\end{frontmatter}

\section{Introduction}
The starting point of this article was a very naive attempt to introduce some geometry, via singularities of differential $1$-forms, in the

\begin{description}
\item[Jacobian problem.]  {\it  Let $F=(P,Q):\mathbb C ^2 \hookleftarrow$ be a polynomial map such that $\frac{\partial (P,Q)}{\partial (x,y)}\neq 0$ on ${\mathbb C}^2$. Then, $F$ is an injective map.}  (See, e.g., \cite{bcw}.)
\end{description}

The relation with differential  $1$-forms is attained by associating to $F$ the differential form $\omega =PdQ-QdP$.

Given $z\in Sing(\omega )$ we have $\omega (z) \wedge dP(z)=\omega (z) \wedge dQ(z)=0$ and, since $dP\wedge dQ = \frac{\partial (P,Q)}{\partial (x,y)} dx\wedge dy \neq 0$, we conclude that $P(z)=Q(z)=0,\;i.e.,\; Sing(\omega )\subset Z(P,Q)$.

On the other hand, since $Z(P,Q) \subset  Sing (\omega )$ trivially holds, we have that $Sing (\omega)=Z(P,Q)$.

This leads at once to the following alternative statement of the Jacobian problem,

\begin{description}
\item[]\hspace{0.25cm} {\it Let $\omega =PdQ-QdP$, where $P,Q$ are polynomials on $\C ^2$. If $d\omega$ has no singular points then, $\omega$ has at most one singular point.}
\end{description}

Maybe that has led us naturally to make use of differential $1$-forms in order to study the

\begin{description}
\item[Jacobian equation.] $\frac{\partial (P,Q)}{\partial (x,y)} =c \in {\mathbb C}$, for $P,Q\in {\mathbb C}[x,y]$.
\end{description}

Which, in fact, has shown to be very efficient in establishing Theorem~\ref{principal}, where a system of algebraic equations involving the homogenoeus components of $P$ and $Q$ is shown to be equivalent to the Jacobian equation.

\section{Preliminaries}
\label{preliminaries}

Henceforward we shall concentrate in investigating the solutions of the Jacobian equation, where $P$ and $Q$ have fixed degrees $k$ and $l$ respectively. In addition, by obvious reasons, $P$ and $Q$ are supposed to satisfy:

\begin{description}

\item[($\imath$)]$P(0)=Q(0)=0$.

\item[($\imath \imath$)] $P,\;Q\neq0$.

\item[($\imath \imath \imath)$]$P$ and $Q$ are not both linear.

\end{description}

Let us now consider the decomposition of $P$ and $Q$ into their respective homogenous components,

$$\begin{array}{c}
P=P_1 +\ldots +P_k\\ \vspace{-.3cm} \\
Q=Q_1 +\ldots +Q_l
\end{array}$$

If $(dP\wedge dQ)_{\mu}$ denotes the homogeneous component of $dP\wedge dQ$ of degree $\mu $, the condition $\frac{\partial (P,Q)}{\partial (x,y)} \in {\mathbb C}$ is equivalent to,

$$(dP \wedge dQ)_{\mu}=0,\;\mu =(k+l)-2, \ldots ,1.$$

Which is, by its turn, equivalent to the following system of $k+l-2$ partial differential equations,

\[ \hspace{0.25cm} \left \{  \begin{array}{cccc}
dP_k \wedge dQ_l=0 \\
dP_k \wedge dQ_{l-1} +dP_{k-1} \wedge dQ_l=0 \\

\vdots \\

dP_2 \wedge dQ_1 +dP_1 \wedge dQ_2=0  \end{array}  \hspace{6cm}    \right . \]
\vspace{0.10cm}

\begin{obs} The above system may be written more conveniently as,
\label{system}
$$({\bf S})\hspace{0.5cm}dP_k \wedge dQ_{l-j} +dP_{k-1} \wedge dQ_{l-(j-1)} + \ldots +dP_{k-j} \wedge dQ_l =0,\; j=0,\ldots k+l-3. $$

Where it is agreed that $P_i=Q_i=0$, whenever $i<0$.

Notice that the $j$-th equation of $({\bf{S}})$ is,
$$(j) \hspace{1cm} \sum \limits_{j\prime =0} \limits^{j}dP_{k-(j-j\prime)} \wedge dQ_{l-j\prime }=0.$$
\end{obs}

Before we proceed to the investigation of the solutions of ({\bf S}), we present below,

\subsection{ Some basic elementary results}
\label{resultadosbasicos}

In what follows, ${\mathbb C}(z)={\mathbb C}(z_1,\ldots ,z_n)$ denotes the field of rational functions on ${\mathbb C}^n$. We shall agree that the zero polynomial is homogeneous of any degree.

\begin{itemize}

\item[(1)]Given a non constant $R\in {\mathbb C}(z)$ we shall denote by $s(R)=max\{m\in \N \mid R=X^m$ for some $X\in \C (z) \}~.$ The notation $G=\sqrt[s]{R}$ means that $s=s(R)$ and that $G^s=R$. Note that necessarily $s(G)=1$, which is equivalent to saying that G is not the power of another rational function. Such a $G$ will be referred to as being {\bf {simple}}.

\item[(2)]
Let $H$ be a holomorphic homogeneous function of degree $k\in \Z$ (defined in some region of ${\mathbb C}^n$). Then, $i ({\mathscr{R}})dH=kH$, where ${\mathscr R}$ denotes the radial vector field  on ${\mathbb C}^n$, i.e., ${\mathscr R}(z)=z,\; z\in {\mathbb C}^n$, and $i({\mathscr R})dH$ is the interior product (see, e.g., \cite{g}, p. 25) of the vector field ${\mathscr R}$ and the differential $1$-form $dH$.

This is just a restatement of the classical Euler's Formula for Homogeneous Functions, in the context of vector fields and differential forms.

\item[(3)]Let $H,J$ be homogeneous holomorphic functions of integer degrees $k,l$ respectively (defined in some region of ${\mathbb C}^n$), such that $H\neq 0$. Then, $dH \wedge dJ=0$ if and only if there exists $\lambda \in \mathbb C$ such that $J^k=\lambda H^l$.

The necessity is an immediate consequence of ($2$) above. Indeed, from the equation $dH\wedge dJ=0$ we have,
$$0=i ({\mathscr R})0=i ({\mathscr R})(dH \wedge dJ)=(i ({\mathscr R})dH)dJ-dH(i({\mathscr R})dJ)=kHdJ-lJdH~.$$

Now, let $M=\frac{J^k}{H^l}$ then $dM=\frac{1}{H^{2l}} (H^l dJ^k -J^k dH^l)=\frac{H^{l-1}J^{k-1}}{H^{2l}}(kHdJ -lJdH)=0$.

Hence, there exists $\lambda \in \mathbb C$ such that $M=\lambda$, i.e., $J^k=\lambda H^l$.

The converse is obvious.

\item[(4)]Let $H,J\in {\mathbb C}(z)$, $H\neq 0$, be quocients of homogeneous polynomials. Then, $dH \wedge dJ=0$ if and only if there exist $\lambda \in \mathbb{C}$ and $t \in \mathbb{Z}$ such that $J=\lambda G^{t}$, where $G=\sqrt[s]{H}$.

In fact, from $G^{s}=H$, we conclude that $G$ is, as well, a quocient of homogeneous polynomials.

If $J=0$ there is nothing to prove. Otherwise, from
$0=dH\wedge dJ=sG^{s-1}dG\wedge dJ$ we deduce that $dG\wedge dJ=0$. By ($3$) we have  $J^g=c G^l$, where $c \in \mathbb C$
and $g,\,l$ are the degrees of $G$ and $J$, respectively.

Now, the result follows by considering the factorizations of the rational functions $G$ and $J$, into irreducible factors, and by noting that $G$ is simple, exactly when the {\it{gcd}}  of the exponents of the factors in its decomposition is equal to $1$.

The converse is evident.

\end{itemize}

\section{Definitions and notation}
\label{resultado}

We shall denote by $\Gamma$ the set of all sequences $\alpha =(\alpha _1,\alpha _2, \ldots)$ of non negative integers, having a finite number of nonzero terms.

Unless otherwise explicitly stated, any sequence appearing in the sequel lies in $\Gamma$.

\newpage
For $\alpha \in \Gamma$ we define,
\begin{itemize}
\item[] \label{modulo} $|\alpha|=\sum \limits _{i\in \N} \alpha _i \,$,
\item[] \label{size}  $\sigma (\alpha )= \sum \limits _{i\in N}i\alpha _i \,$,
\item[] \label{dominio} $D_{\alpha}=\left\{i\in \N \mid \alpha _i \neq 0 \right\} \,$.
\end{itemize}

The functions $|\alpha|$ and $\sigma (\alpha)$ will be referred to, respectively, as the modulus and the size of $\alpha $.\\

For each $j\in \N$, we denote by $e_j$ the sequence whose $j$-th term is $1$ and all the others are zero.

Given a nonzero sequence $\alpha$, for each $i\in D_{\alpha}$, we denote by $\alpha (i)$ the sequence $\alpha - e_i$. The function $i\in D_{\alpha}\longmapsto \alpha (i)\in \Gamma$ will be referred to as the function $\alpha (i)$.\\

If $t\in \R$ and $\alpha=(\alpha_1,\alpha _2,\ldots ,\alpha _k,0,0,\ldots)\in \Gamma$, we shall denote by ${t \choose  \alpha}$ the usual  multinomial coefficient ${t \choose  \alpha_1,\ldots ,\alpha _k}=\frac{(t)_{|\alpha|}}{\alpha !}$, where $\alpha != \prod \limits _{i\in \mathbb N}\alpha _i !$ and $(t)_{|\alpha|}$ is the Pochhamer symbol for the falling factorial $t(t-1)\ldots (t-|\alpha|+1)$. Recall that when $\alpha =0,\;(t)_0 =1$ by definition.

We point out, for further reference, the following elementary,

\begin{id}
\label{id}
 \Large{${t \choose  \alpha + e_i}=\frac{(t-|\alpha|)}{\alpha _i +1}{t \choose  \alpha}$}.
\end{id}

Finally, let $X=(X_1,X_2,\ldots)$, where $X_1,X_2,\dots$ are indeterminates. Given $\alpha \in \Gamma $ we shall adopt the usual notation $X^{\alpha}=X_1 ^{\alpha _1}X_2 ^{\alpha _2}\ldots$.

\section{Statement and proof of the result}

\begin{teo}
\label{principal}
Let $P, Q\in \C [x,y]$, of degrees $k$ and $l$ respectively,  be such that $P(0)=Q(0)=0$, and $kl>1$. Then, $\frac{\partial (P,Q)}{\partial (x,y)} \in {\mathbb C}$ if and only if there exist unique  $\lambda _r \in \C,\; r=0,\ldots ,k+l-3$, such that,
$$Q_{l-j} = \sum \limits _{r=0} \limits ^{j}\sum \limits _{\sigma(\alpha) = j - r}\lambda _r {{s_r /s}\choose  \alpha} G ^{s_r - s|\alpha|} P\_ ^{\alpha}, \;\; 0 \leq j \leq k+l-3, $$
\noindent where $G=\sqrt[s]{P_k};\;P\_ = (P_{k-1}, P_{k-2} , .... )$;  $s_r = s\frac{(l - r)}{k}$, if  $\lambda_r \neq 0$,  and $s_r=0$, if  $\lambda_r = 0$. Furthermore,  $s_r$ turns out to be an integer, whenever  $\lambda _r \neq 0$.
\end{teo}

\noindent{\bf Proof.}

We shall omit, along the proof, details that turn out to be mere elementary algebraic manipulations.

Henceforth, in order to simplify the typing, and the reading, we set $t_r =s_r /s$.

We shall see that the above expressions of $Q_{l-j}$ are obtained by solving recursively all equations of the system {\bf (S)} for the $Q_{l-j}$.

In fact, we will show, by recurrence on $j$, the following assertion:

\begin{description}
\item[({A$_{\bf j}$})] Given $0\leq j\leq k+l-3$, the first equations of the system ({\bf{S}}) up to the $j$- th, hold iff there exist unique $\lambda _r \in \C,\;r=0,\ldots ,j,$ such that,

\hspace{1.5cm}$Q_{l-j\prime } = \sum \limits _{r=0} \limits ^{j\prime}\sum \limits _{\sigma(\alpha) = j\prime - r}\lambda _r {{t_r }\choose  \alpha} G ^{s_r - s|\alpha|} P\_ ^{\alpha}, \;\;  j\prime =0,\ldots ,j.$
\end{description}

For $j=0$, the assertion is a straightforward consequence of ($4$) in subsection~\ref{resultadosbasicos}.

In order to complete the recurrence procedure let us prove that {\bf (A$_{\bf {j-1}}$)} implies
{(\bf A$_{\bf j}$)}, for $0<j\leq k+l-3$.

Indeed, by the recurrence hypothesis, {(\bf A$_{\bf j}$)} is equivalent to:\\

The $j$-th equation of the system holds iff there exists a unique $\lambda _j \in \C$ such that,
$$
Q_{l-j} = \sum \limits _{r=0} \limits ^{j}\sum \limits _{\sigma(\alpha) = j - r}\lambda _r {{t_r}\choose  \alpha} G ^{s_r - s|\alpha|} P\_ ^{\alpha}\;\;,
$$
where $s_j = s\frac{(l - j)}{k}$, if  $\lambda_j \neq 0$,  and $s_j=0$, if  $\lambda_j = 0$. Moreover,  $s_j \in \Z$, if  $\lambda _j \neq 0$.\\

By Remark~\ref{system}, the $j$-th equation of ({\bf S}) is,
$$(j)\hspace{1.5cm}dP_k \wedge dQ_{l-j} + \sum \limits_{j\prime =0} \limits^{j-1}dP_{k-(j-j\prime)} \wedge dQ_{l-j\prime }=0.$$

Since $P_k=G^s$, we have that $dP_k\wedge dQ_{l-j}=sG^{s-1}dG \wedge dQ_{l-j}$.\\

Now, we shall compute $dQ_{l-j\prime}$ and $dP_{k-(j-j\prime)} \wedge dQ_{l-j\prime}$.\\

By the recurrence hypothesis we have,\\ \vspace{.2cm}
$dQ_{l-j\prime }=d[  \sum \limits _{r=0} \limits ^{j\prime}\sum \limits _{\sigma(\alpha) = j\prime - r}\lambda _r {{t_r}\choose  \alpha} G ^{s_r - s|\alpha|} P\_ ^{\alpha}  ]=\\
= \sum \limits _{r=0} \limits ^{j\prime}\sum \limits _{\sigma(\alpha) = j\prime - r}\lambda _r {{t_r}\choose  \alpha}[ (s_r -s|\alpha|)G ^{s_r - s|\alpha|-1} P\_ ^{\alpha}dG+ G ^{s_r - s|\alpha|} dP\_ ^{\alpha}  ]$.\\ \vspace{.2cm}

By taking the exterior product of $dP_{k-(j-j\prime)}$ and the expression above, we obtain,\\
$dP_{k-(j-j\prime)}\wedge dQ_{l-j\prime}=$\\
$\sum \limits_{r=0} \limits^{j\prime}
\sum \limits_{\sigma(\alpha) = j\prime - r}
\lambda _r {{t_r}\choose  \alpha}
\left[(s_r -s |\alpha|)G ^{s_r - s|\alpha|-1} P\_ ^{\alpha}dP_{k-(j-j\prime)}\wedge dG~+~G ^{s_r - s|\alpha|} dP_{k-(j-j\prime)}\wedge dP\_ ^{\alpha}\right],$ \vspace{.1cm}\\
$0\leq j\prime \leq j-1$.\\

Taking into account the above expressions, equation {\it (j)} is now,\\
$sG^{s-1}dG \wedge dQ_{l-j}+\\
\sum \limits _{j\prime =0} \limits ^{j-1}\sum \limits _{r=0} \limits ^{j\prime}\sum \limits _{\sigma(\alpha) = j\prime - r}\lambda _r {{t_r}\choose  \alpha}(s_r -s|\alpha|)G ^{s_r - s|\alpha|-1} P\_ ^{\alpha}dP_{k-(j-j\prime)}\wedge dG+\\
 \sum \limits _{j\prime =0} \limits ^{j-1}\sum \limits _{r=0} \limits ^{j\prime}\sum \limits _{\sigma(\alpha) = j\prime - r}\lambda _r {{t_r}\choose  \alpha}G ^{s_r - s|\alpha|} dP_{k-(j-j\prime)}\wedge dP\_ ^{\alpha} =0$.\\

By factoring out $sG^{s-1}dG$, we obtain,\\
$sG^{s-1}dG \wedge [dQ_{l-j}-\sum \limits _{j\prime =0} \limits ^{j-1}\sum \limits _{r=0} \limits ^{j\prime}\sum \limits _{\sigma(\alpha) = j\prime - r}\lambda _r {{t_r}\choose  \alpha}(t_r - |\alpha|)G ^{s_r - s(|\alpha| +1)} P\_ ^{\alpha}dP_{k-(j-j\prime)} ]+\\
 \sum \limits _{j\prime =0} \limits ^{j-1}\sum \limits _{r=0} \limits ^{j\prime}\sum \limits _{\sigma(\alpha) = j\prime - r}\lambda _r {{t_r}\choose  \alpha}G ^{s_r - s|\alpha|} dP_{k-(j-j\prime)}\wedge dP\_ ^{\alpha} =0$.\\

Now, the recurrence procedure follows easily from the two statements below,

\begin{statement}
\label{statement1}
$\vspace{.2cm} \\ \sum \limits _{j\prime =0} \limits ^{j-1}\sum \limits _{r=0} \limits ^{j\prime}\sum \limits _{\sigma(\alpha) = j\prime - r}\lambda _r {{t_r}\choose  \alpha}G ^{s_r - s|\alpha|} dP_{k-(j-j\prime)}\wedge dP\_ ^{\alpha} =0,\;0 < j \leq k+l- 3$.
\end{statement}

\begin{statement}
\label{statement2}
$\vspace{.2cm} \\ \sum \limits _{j\prime =0} \limits ^{j-1}\sum \limits _{r=0} \limits ^{j\prime}\sum \limits _{\sigma(\alpha) = j\prime - r}\lambda _r {{t_r}\choose  \alpha}(t_r -|\alpha|)G ^{s_r - s(|\alpha|+1)} P\_ ^{\alpha}dP_{k-(j-j\prime)}=
\sum \limits _{r=0} \limits ^{j-1}\sum \limits _{\sigma(\alpha) = j - r}\lambda _r {{t_r}\choose  \alpha} G ^{s_r - s|\alpha|} dP\_ ^{\alpha}$,\\

\noindent $\;0 < j \leq k+l - 3$.
\end{statement}

As a matter of fact, by taking for granted the two statements above, the equation {\it (j)} becomes,\\
$0=dG \wedge [ dQ_{l-j}-\sum \limits _{r=0} \limits ^{j-1}\sum \limits _{\sigma(\alpha) = j - r}\lambda _r {{t_r}\choose  \alpha} G ^{s_r - s|\alpha|} dP\_ ^{\alpha} ]=\\
=dG \wedge d[ Q_{l-j}-\sum \limits _{r=0} \limits ^{j-1}\sum \limits _{\sigma(\alpha) = j - r}\lambda _r {{t_r}\choose  \alpha} G ^{s_r - s|\alpha|} P\_ ^{\alpha} ]$, once $dG\wedge dG=0$.\\

We notice that $deg(G ^{s_r - s|\alpha|} P\_ ^{\alpha}) =l-j$, if $\sigma (\alpha )=j-r$ and $\lambda _r \neq 0$. Thus, $Q_{l-j}-\sum \limits _{r=0} \limits ^{j-1}\sum \limits _{\sigma(\alpha) = j - r}\lambda _r {{t_r}\choose  \alpha} G ^{s_r - s|\alpha|} P\_ ^{\alpha}$ is a quocient of homogeneous polynomials.

Now, since $G$ is simple, it follows from ($4$) of subsection~\ref{resultadosbasicos}, that the equation $(j)$ holds iff,
\begin{center}
$Q_{l-j}-\sum \limits _{r=0} \limits ^{j-1}\sum \limits _{\sigma(\alpha) = j - r}\lambda _r {{t_r}\choose  \alpha} G ^{s_r - s|\alpha|} P\_ ^{\alpha}=\lambda _j G^{s_j}$, for some $\lambda _j \in \C$ and $s_j \in \Z$.
\end{center}

Clearly, the constant $\lambda _j$ is uniquely determined by the above equation and,
if $\lambda _j \neq 0$, this equation implies that $l-j=s_jdeg(G)=s_j \frac{k}{s}$, i.e., $s_j=\frac{s}{k}(l-j)$.

On the other hand, if $\lambda _j=0$, we may clearly choose $s_j=0$.

In other words, we have just shown that, under the recurrence hypothesis, the identity,

\begin{center}
$
Q_{l-j} = \sum \limits _{r=0} \limits ^{j-1}\sum \limits _{\sigma(\alpha) = j - r}\lambda _r {{t_r}\choose  \alpha} G ^{s_r - s|\alpha|} P\_ ^{\alpha}+\lambda _j G^{s_j} = \sum \limits _{r=0} \limits ^{j}\sum \limits _{\sigma(\alpha) = j - r}\lambda _r {{t_r}\choose  \alpha} G ^{s_r - s| \alpha|} P\_ ^{\alpha}\; ,
$
\end{center}

\noindent with $\lambda _j$ and $s_j$ as described above, is in fact equivalent to equation $(j)$.

\qed

Now, we will provide the proof of the two statements.

\subsection{Proof of Statement~\ref{statement1}}

Before we proceed to the proof we set,
\begin{center}
$A= \sum \limits _{j\prime =0} \limits ^{j-1}\sum \limits _{r=0} \limits ^{j\prime}\sum \limits _{\sigma(\alpha) = j\prime - r}\lambda _r {{t_r}\choose  \alpha}G ^{s_r - s|\alpha|} dP_{k-(j-j\prime)}\wedge dP\_ ^{\alpha}$.
\end{center}

We are supposed to prove that  $A = 0$.

\newpage
First we remark that $dP\_ ^{\alpha} = 0$, if $\alpha = 0$ and, otherwise,  $dP\_ ^{\alpha} = \sum \limits _{i \in D _{\alpha}} \alpha _i P\_ ^{\alpha (i)}dP_{k-i}$.\\

Consequently, we have that $dP_{k-(j-j')} \wedge dP\_ ^{\alpha} = 0$, if  $\alpha = 0$  and,   $dP _{k-(j-j')}\wedge dP\_ ^{\alpha}  = \sum \limits _{i \in D _{\alpha}} \alpha _i P\_ ^{\alpha (i)}dP _{k-(j-j')}\wedge dP_{k-i}$,  if  $\alpha\neq 0$.\\

Hence, in the expression of  $A$,  we may restrict ourselves to those summands where  $\alpha \neq 0$, if any exist. If not, $A$ trivially vanishes.

The condition  $\alpha \neq 0$  may be more appropriately expressed in terms of the indexes range, by observing that, $\alpha \neq 0  	 \Longleftrightarrow \sigma(\alpha) \neq 0  	\Longleftrightarrow j' - r \neq 0  	\Longleftrightarrow  0 \leq r \leq j' - 1 \leq j - 2$,  which is equivalent to  $1 \leq j' \leq j - 1,\;  0\leq r \leq j' - 1$  and  $j\geq 2$.\\

The above discussion may be summarized as follows:\\

$A = 0$, if  $j < 2$, and for  $j \geq 2$ we have,
\begin{center}
$(\star) \hspace{.5cm}A= \sum \limits _{j\prime =1} \limits ^{j-1}\sum \limits _{r=0} \limits ^{j\prime -1}\sum \limits _{\sigma(\alpha) = j\prime - r}\sum \limits _{i\in D_{\alpha}} \lambda _r {{t_r}\choose  \alpha}G ^{s_r - s|\alpha|}\alpha _i P\_ ^{\alpha (i)} dP_{k-(j-j\prime)}\wedge dP_{k-i}$.
\end{center}

Henceforward we shall presume  $j \geq 2$.

Let us denote by  $\mathbb A$  the set of all  $4$-tuples  $a=\begin{pmatrix}\,j\prime \,\\ r \\ \alpha \\ i\end{pmatrix}$, whose coordinates are subjected to the same constraints specified in ($\star$) above.

Clearly,\\

$A=\sum \limits _{a\in \mathbb A} \Phi (a)$, where $\Phi (a)=\lambda _r {{t_r}\choose  \alpha}G ^{s_r - s|\alpha|}\alpha _i P\_ ^{\alpha (i)} dP_{k-(j-j\prime)}\wedge dP_{k-i}$.

Now we set, for  $a \in \mathbb A,\; \tau(a) = \begin{pmatrix}\,j-i \,\\ r \\ \alpha (i) + e_{j-j\prime} \\ j-j\prime\end{pmatrix}$. It is immediate to check that this defines, in fact, a bijective function $\tau:\mathbb A\longrightarrow \mathbb A$. Such function  satisfies:
$$\Phi (\tau (a))=-\Phi (a)~.$$

\newpage
As a matter of fact, by the very definitions of   $\Phi$  and  $\tau$ we have, $\Phi(\tau (a))=$\vspace{.2cm} \\
$\lambda _r {{t_r}\choose  \alpha (i)+e_{j-j\prime}}G ^{s_r - s|\alpha (i)+e_{j-j\prime}|}(\alpha (i)+e_{j-j\prime})_{j-j\prime} P\_ ^{(\alpha (i)+e_{j-j\prime})(j-j\prime)} dP_{k-i}\wedge dP_{k-(j-j\prime)},$\vspace{.2cm} \\
and then, the fact that  $\Phi (\tau(a)) = -\Phi (a)$  turns out to be an immediate consequence of Identity~\ref{id}.\\

Hence we conclude that
$$A=\sum \limits_{a \in \mathbb{A}} \Phi(a)= \sum \limits_{a \in \mathbb{A}} \Phi(\tau(a))=-A~, ~\mbox{i.e.}~A=0~.$$
\qed

\subsection{Proof of Statement~\ref{statement2}}
Let us set,

$B=\sum \limits _{j\prime =0} \limits ^{j-1}\sum \limits _{r=0} \limits ^{j\prime}\sum \limits _{\sigma(\alpha) = j\prime - r}\lambda _r {{t_r}\choose  \alpha}(t_r -|\alpha|)G ^{s_r - s(|\alpha| +1)} P\_ ^{\alpha}dP_{k-(j-j\prime)}$,

$B'=\sum \limits _{r=0} \limits ^{j-1}\sum \limits _{\sigma(\alpha) = j - r}\lambda _r {{t_r}\choose  \alpha} G ^{s_r - s|\alpha|} dP\_ ^{\alpha}=\sum \limits _{r=0} \limits ^{j-1}\sum \limits _{\sigma(\alpha) = j - r}
\sum \limits _{i\in D_{\alpha}}
\lambda _r {{t_r}\choose  \alpha} G ^{s_r - s|\alpha|} \alpha _iP\_ ^{\alpha (i)}dP_{k-i}$,
this last equality is due to the fact that  $\alpha \neq 0$, once  $\sigma (\alpha) = j - r \neq 0$.\\

Recall we want to prove that  $B = B'$.

The proof consists basically in showing that the summands in the expressions of  $B$  and  $B'$  are exactly the same. To this end, we shall express both, $B$  and  $B'$,  into the more suitable form:

$B=\sum \limits _{b\in \mathbb B}\Psi (b)$, where
$\mathbb B=\left\{\begin{pmatrix}\,j\prime \,\\ r \\ \alpha \end{pmatrix}\mid \; 0\leq j\prime \leq j-1,\;0\leq r\leq j\prime ,\;\sigma (\alpha)=j\prime -r \right\}$, \linebreak and  $\Psi (b)=\lambda _r {{t_r}\choose  \alpha}(t_r -|\alpha|)G ^{s_r - s(|\alpha| +1)} P\_ ^{\alpha}dP_{k-(j-j\prime)}$.

\vspace{0.35cm}

$B'=\sum \limits _{b'\in \mathbb B '}\Psi ' (b')$, where $\mathbb B '= \left\{\begin{pmatrix}\,r  \,\\ \alpha \\ i  \end{pmatrix} \mid  \; 0\leq r\leq j-1 ,\;\sigma (\alpha)=j -r ,\; i\in D_{\alpha } \right\}$, \linebreak and $\Psi ' (b')=\lambda _{r} {{t_{r}}\choose  \alpha }G ^{s_{r} - s|\alpha|} {\alpha }_i P\_ ^{\alpha  (i)}dP_{k-i}$.

Now, for $b\in \mathbb B$, we set  $\varrho (b) =  \begin{pmatrix} \, r\,  \\  \alpha + e_{j-j\prime}   \\ j-j\prime \end{pmatrix}$. It can be easily verified that this defines a bijective function   $\varrho :\mathbb B \longrightarrow \mathbb B '$.

Obviously, in order to conclude the proof, it suffices to show that  $\Psi(b) = \Psi '(\varrho (b))$. Indeed,
$$B' = \sum \limits _{b'\in \mathbb B '}\Psi '(b') = \sum  \limits _{b'\in \varrho (\mathbb B )}\Psi '(b') = \sum  \limits _{b\in \mathbb B }\Psi '(\varrho (b)) = \sum  \limits _{b\in \mathbb B }\Psi (b) = B~.$$

Finally, by a direct computation we find that,
$$\Psi '(\varrho (b)) = \lambda _{r} {{t_{r}}\choose  \alpha + e_{j-j\prime}}G ^{s_{r} - s|\alpha + e_{j-j\prime}|} {(\alpha + e_{j-j\prime}) }_{j-j\prime} P\_ ^{(\alpha + e_{j-j\prime}) (j-j\prime)}dP_{k-(j-j\prime)}$$
\noindent And then, the fact that  $\Psi '(\varrho (b))=\Psi(b)$ follows at once from Identity~\ref{id}.\\
\qed \vspace{.3cm}

\section{Final comments}
It is worth mentioning that, when  $k > 1$, the number of equations, provided by Theorem~\ref{principal}, is  $k + l - 2  \geq  l$. Hence, the first  $l$  equations are explicit expressions of the homogeneous components of  $Q$  in terms of those of  $P$, whereas the remaining equations involve only the homogeneous components of  $P$.

In particular, this holds when  $k \geq l$, and, as far as our purpose is concerned, we could have restricted ourselves to this case, by simply reordering the pair  $(P,Q)$, if necessary.

We have made our choice for the current statement of the theorem, mostly because of its "symmetric" character:

When  $P$  and  $Q$  are not linear we can, indistinctly, express the homogeneous components of  $Q$  in terms of those of  $P$, and conversely.

We point out that the theorem clearly holds in dimension $n \geq 2$, by replacing  {$\frac{\partial (P,Q)}{\partial (x,y)} \in {\mathbb C}$  by,  $dP\wedge dQ $ is constant. The particular case  $dP\wedge dQ = 0$  corresponds to the problem of algebraic dependence of the polynomials  $P,\; Q$ (see, e.g., \cite{pi}, Lemma 1, and \cite{hp}, Ch. III).

Finally, it is evident that, mutatis mutandis, a real version of  Theorem 4.1  is promptly available. For  $c \neq 0$, it is related to the real Jacobian conjecture (see, e.g., \cite{bcw}, Part II, 10.1, and \cite{sergey}).

\newpage
\section*{List of notations}
$\C$, the field of complex
numbers.

$\N=\left\{1,2,\ldots \right\}$, the set of natural numbers.

$\mathbb Z$, the set of integer numbers.

${\mathbb C}[x,y]$, the ring of polynomials on ${\mathbb C}^2$.

$Sing(\omega)=\left\{z\in \C ^2\mid \omega (z)=0 \right\}$, the set of singularities of the differential form $\omega$.

$Z(P,Q)=\left\{z\in \C ^2\mid P(z)=Q(z)=0 \right\}$, the set of zeros of the mapping $(P,Q)$.

$deg(H)$, the degree of the homogeneous function $H$.

\section*{References}

\stepcounter{section}

\end{document}